\documentclass[a4, 11.5pt]{article}
\usepackage{amssymb,latexsym,amsmath, mathtools, array, graphicx, authblk, enumerate}
\usepackage{amsmath, amssymb , graphicx, fancyvrb, verbatimbox,placeins}
\usepackage[shortlabels]{enumitem}
\usepackage[font=singlespacing]{caption}
\usepackage{xcolor}
\usepackage{geometry,graphicx}
\newtheorem{definition}{Definition}[section]
\numberwithin{equation}{section}
\newcolumntype{L}{>{\centering\arraybackslash}m{3cm}}
\begin{document}
	\date{}
\title{\bf On Solving Fixed Charge Transportation Problems Having Interval Valued Parameters}
\author{ Ummey Habiba$^{1}$\thanks{Email : say2habiba@gmail.com},
	~
	Abdul Quddoos$^{2}$\thanks{Email : abdulquddoos.or.amu@gmail.com}
	~and	Masihuddin$^{3}$\thanks{Email : masih.iitk@gmail.com, masihuddin@iiitg.ac.in}}
\maketitle \noindent {\it $^{1}$ Department of Mathematics, Integral University, Lucknow-226026
	, Uttar Pradesh, India} 
\\ {\it $^{2}$ Department of Mathematics, Era University, Lucknow-226003, Uttar Pradesh, India} 
\\  {\it $^{3}$ Department of  Science \& Mathematics, Indian Institute of Information Technology Guwahati-781015, India
}
\vspace{4mm}\\
%\newcommand{\oddhead}{ Estimation of the Selected Treatment Mean}
%\renewcommand{\@oddhead}
%{\hspace*{-3pt}\raisebox{-3pt}[\headheight][0pt]
%	{\vbox{\hbox to \textwidth
		%	{\hfill\oddhead}\vskip8pt}}}
%\vspace*{0.05in}
\noindent {\bf Abstract}:
%\title{On solving an interval fixed-charge transportation problem}
%\author{Ummey Habiba, Abdul Quddoos}
%
%\author{Masihuddin}
%\affil {Integral University, Lucknow, India}
%\maketitle
In this article, we propose a new method for solving the interval fixed charge transportation problem (IFCTP), wherein the parameters (associated cost, fixed cost, supply, and demand) are represented by interval numbers. First, an equivalent bi-objective fixed charge transportation problem (FCTP) is derived from the given IFCTP, and then the equivalent crisp problem is solved using a fuzzy programming technique. To demonstrate the solution procedure, two existing numerical examples (Safi and Razmjoo {\cite{bakp1}}) are coded and solved in LINGO 19.0. We establish the effectiveness of our proposed method through a comparison of the results achieved with those of two pre-existing methods.
\vspace{2mm}\\
\noindent {\it AMS 2010 SUBJECT CLASSIFICATIONS:} 90C08 · 90C70   \\

\noindent {\it Keywords and Phrases:}~Fixed Charge Transportation Problem; Interval Transportation Problem; Fuzzy Programming.
%\end{abstract}
\let\thefootnote\relax\footnote{Corresponding Author Email: masih.iitk@gmail.com }
 %{\it{Key words}} : Fixed charge Transportation Problem; Interval Transportation Problem; Fuzzy Programming
\section{Introduction}
 In a traditional transportation problem, the objective is to determine the least expensive method of transporting a homogeneous product from $m$ origins to $n$ destinations. In reality, however, there may be additional charges due to tax fees, opening routes, property tax, airport landing fees, etc., which may be referred to as a fixed cost (or fixed charge) that is independent of the quantity of product being transported. A new variant of the classic transportation problem known as Fixed Charge Transportation Problem (FCTP) emerges when fixed costs are also considered. In classical transportation problems, the parameters are precisely known to the decision-maker. However, in many practical situations, these parameters may reflect imprecise nature due to lack of evidence, globalization of the market, insufficient information etc. This type of imprecise data can be handled very well by using interval numbers.
The Interval Fixed Charge Transportation Problem (IFCTP) is another variant of the traditional transportation problem in operations research. In the IFCTP, you have a set of suppliers and a set of customers, each with a fixed supply or demand quantity. The objective is to determine how to transport goods from suppliers to customers while minimizing the total cost. The cost of transportation typically includes both a variable cost based on the amount shipped and a fixed cost associated with setting up the transportation.
The unique aspect of the IFCTP is that it deals with uncertainty in the cost parameters. Instead of having exact transportation costs, you have interval numbers representing the possible range of costs for each shipment. This uncertainty can arise due to factors such as fluctuations in fuel prices, market conditions, and other external variables that affect transportation costs.
Solving the IFCTP involves finding an allocation plan that specifies how much to ship from each supplier to each customer while considering the uncertain interval costs. The goal is to find a solution that minimizes the total cost while accounting for the uncertainty in the cost parameters. 
 Initially, Hirsch and Dantzig {\cite{bakp2}} had developed FCTP as a generalization of the conventional transportation problem.\par 
 There are several approaches of solving the IFCTP. Ishibuchi and Tanaka {\cite{bakp3}} were the first who proposed the order relations of interval numbers which was later generalized by Chanas and Kuchta {\cite{bakp4}}. Later on, Sengupta and Pal {\cite{bakp5}} defined an acceptability function for comparison of interval numbers. Hu \textit{et} al. {\cite{bakp6}} redefined the order relations of interval numbers given by Ishibuchi and Tanaka {\cite{bakp3}}. Using the order relations of Ishibuchi and Tanaka {\cite{bakp3}}, Das \textit{et} al. {\cite{bakp7}} developed a method to solve interval transportation problem in which they minimize the center and the right limit of the interval valued objective function simultaneously. Sengupta and Pal {\cite{bakp8}} proposed another fuzzy orientation method for solving interval transportation problem by simultaneously minimizing the center and width of the interval valued objective function. Panda and Das {\cite{bakp9}} proposed an algorithm for solving a two vehicle cost varying transportation problem . Safi and Razmjoo{\cite{bakp1}} worked on an IFCTP using the order relations of (Ishibuchi and Tanaka {\cite{bakp3}} and Hu \textit{et} al. {\cite{bakp6}}). Zhang \textit{et} al.{\cite{bakp11}} investigated fixed charge solid transportation problem in uncertain environment. Roy and Midya {\cite{bakp10}} presented an uncertain intuitionistic fuzzy model for multi-objective fixed charge solid transportation problem. Sifaoui and Aïder {\cite{bakp12}} proposed an interval multi-objective multi-item fixed charge solid transportation problem along with budget constraint and safety measure.\par 
This paper proposes a novel optimistic method for solving IFCTP in which the interval-valued objective function is converted into two crisp objective functions. The first objective function minimizes the left limit (i.e., the best-case scenario), whereas the second objective function minimizes the width of the interval-valued objective function (i.e., uncertainty). Using fuzzy programming, the Pareto optimal solution to the equivalent crisp bi-objective problem is obtained. For the purpose of demonstrating the effectiveness of the proposed method, two numerical examples were considered and solved using LINGO 19.0. The paper contains a detailed solution to one numerical example along with its Lingo codes for convenience of the readers. At the end of the paper, the proposed method has been compared with the existing method of  Safi and Razmjoo {\cite{bakp1}}.\\\par 
In the next section we provide some preliminary definitions that will be used in the later sections of the paper.
\section{Preliminaries}
\begin{definition}[{\bf{Interval Numbers} \cite{bakp5}}]\label{definition1}
Let the small letters $i$, $j$ represent real numbers and the capital letters $I$, $J$  represent the closed intervals on the real line $\mathbb{R}$. The intervals $I$ and $J$ are defined as:\\
\begin{equation*}
I=[\underline{i}, ~\overline{i}]={\{\ i: \underline{i}~ \leq\ i \leq\ \overline{i},~ i\in\mathbb{R} \}}\\
\end{equation*}
and \begin{equation*}
	J=[\underline{j}, ~\overline{j}]={\{\ j: \underline{j}~ \leq\ j \leq\ \overline{j},~ j\in\mathbb{R} \}},\\
\end{equation*}
where $\underline{i}$ ($\underline{j}$) and $\overline{i}$ ($\overline{j}$) are the left-limit and right-limit of an interval $I$ (J) on the real line $\mathbb{R}$.
The interval $I$ can be represented in termS of center and half-width (it is simply called width), as,
\begin{equation*}\label{2}
I=\big \langle{i_{c}, i_{w}}\big \rangle={\{\ i:\ i_{c}-i_{w}\leq\ i\leq\ i_{c}+i_{w}, i\in\mathbb{R} \}},\\
\end{equation*}
where, $i_{c}=\big(\frac{\overline{i}+\underline{i}}{2}\big)$ and $i_{w}=\big(\frac{\overline{i}-\underline{i}}{2}\big)$ are the center and width of an interval $I$. Similar representation can be given for the interval $J$.
\end{definition}
\begin{definition}[{\bf{ Interval Arithmetic }\cite{bakp5}}]\label{definition2}
 Let $I=[\underline{i}, ~\overline{i}]$ and $J=[\underline{j}, ~\overline{j}]$ be the two closed intervals as defined above on the real line $\mathbb{R}$. Then,
\begin{align*}
&I+J=[\underline{i}, ~\overline{i}]+[\underline{j}, ~\overline{j}]=[\underline{i}+\underline{j},~ \overline{i}+\overline{j}]\\
&I+J=\big \langle{i_{c}, ~i_{w}}\big \rangle+\big \langle{j_{c},~ j_{w}}\big \rangle=\big \langle{i_{c}+j_{c}, ~i_{w}+j_{w}}\big \rangle\\
&{\gamma}I=\gamma[\underline{i}, ~\overline{i}]=
\begin{dcases}
[{\gamma}\underline{i},~ {\gamma}\overline{i}]~&\text{if}~\gamma{\geq}0\\
[{\gamma}\overline{i},~ {\gamma}\underline{i}]~&\text{if}~\gamma{<}0
\end{dcases}\\
&{\gamma}I=\gamma{\big \langle{i_{c}, {i_{w}}}\big \rangle}={\big \langle{\gamma}{i_{c}, {|\gamma|}{i_{w}}}\big \rangle}
\end{align*}
where, $\gamma$ is a real number.\\
\end{definition}
\subsection{A new order relation for comparing two intervals}
In this section, we propose the concept of order relation between two intervals, which will be utilized later in this paper to compare the cost of two interval fixed-cost transportation problems. The preference among two intervals can be decided with the help of an order relation $``\leq_{D}"$ which is defined as follows:\par
\noindent  
\textbf{Definition}: \label{def1} Let $P$ and $Q$ be two intervals which represent uncertain costs from two alternatives. Consider the cost of each alternative lie in the corresponding interval.\\
The order relation $\leq{_{D}}$ between $P={\big \langle{p_{c}, {p_{w}}}\big \rangle}$ and $Q={\langle{q_{c}, q_{w}}\big \rangle}$ is defined as:
\begin{equation*}
	%\left\{
	\begin{array}{@{}ll@{}}
		P\leq_{D}{Q}~~\text{if}~~d_{IP}\leq d_{IQ}\\
		P<_{D}{Q}~~\text{if}~~P\leq_{D}{Q}~~\text{and}~~ P\neq Q
	\end{array}%\right.
\end{equation*}
where, $I={\langle{i_{c},i_{w}}\big \rangle}$ represent the ideal expected-value and ideal uncertainty.
\begin{equation*}
	\begin{array}{@{}ll@{}}
		d_{IP}=\sqrt{(p_{c}-i_{c})^{2}+(p_{w}-i_{w})^{2}} \\
		d_{IQ}=\sqrt{(q_{c}-i_{c})^{2}+(q_{w}-i_{w})^{2}}
	\end{array}
\end{equation*}
If $P\leq_{D}{Q}$, then $P$ is preferred over $Q$.\vspace{3mm}

In the next section we provide a detailed algorithm for solving the interval fixed charge transportation problem (IFCTP). 
\section{An alternative approach of solving interval-valued transportation problems}

 \subsection{Mathematical models of FCTP and IFCTP}
 The general mathematical models of FCTP and IFCTP has been discussed in the following subsections (\ref{4.1}) and (\ref{4.2}), respectively.
 \subsubsection{Mathematical model of FCTP} \label{4.1}
The FCTP can be seen as a transportation problem with $m$ sources (factories) and $n$ destinations (customers). The product can be shipped from each source to any of the destinations for $t_{ij}$ per unit plus a fixed cost of $l_{ij}$. In a balanced FCTP, the total amount of product available at all sources must be equal to the total amount of product needed at all destinations. But in real life, these conditions may not always be met. For this reason, we have employed the standard form of FCTP instead of the balance form, which can be written mathematically as:\\ 

{\bf{Problem-I:}}
 \begin{align*}
   &{\rm Minimize~} Z=\sum_{i=1}^{m} \sum_{j=1}^{n}({t}_{ij}y_{ij}+{l}_{ij}x_{ij}) \\
   {\rm subject~to~the~constraints;}\nonumber& \\
   &\sum_{j=1}^{n}y_{ij} \leq {\alpha}_{i}, ~i=1,2,\ldots,{m} \\
   &\sum_{i=1}^{m}y_{ij} \geq {\beta}_{j}, ~j=1,2,\ldots,{n}\\
   &y_{ij}\geq{0},~i=1,2,\ldots,{m},~ j=1,2,\ldots,{n}\\
   &x_{ij} =
\begin{dcases}
       0~&{\rm if}~y_{ij}=0\\
       1~&{\rm if}~y_{ij}>0
\end{dcases},\\
{\rm{model~ is~ feasible~ if}},&\sum_{i=1}^{m}{\alpha}_{i} \geq \sum_{j=1}^{n}{\beta}_{j}
   \end{align*}
   where,
   \begin{align*}
  {t}_{ij} &: \text{transportation cost of unit product from the $i^{th}$ source  to the $j^{th}$ destination} \\
  {l}_{ij} &: \text{fixed charge which occurs between the the $i^{th}$ source and the  $j^{th}$ destination} \\
 {\alpha}_{i} &: \text{availability at the $i^{th}$ source}\\
 {\beta}_{j}&: \text{demand at the  $j^{th}$ destination}\\
 y_{ij}&: \text{quantity transported from the $i^{th}$ source  to the $j^{th}$ destination}\\
 x_{ij}&: \text{binary variable}.
 \end{align*}
\subsubsection{Mathematical model of IFCTP}\label{4.2}
 In Problem-I it is assumed that all the parameters ${t}_{ij}$, ${l}_{ij}$, ${\alpha}_{i}$, ${\beta}_{j}$ are known crisp quantities. Still, in practical situations, one may not know about the crisp values of these parameters, but a range in which the values of these parameters lie can be estimated based on past data. So, if we replace these crisp parameters in Problem-I with an interval number, it converts into IFCTP. The mathematical representation of IFCTP can be given as follows:\\ 
{\bf{Problem-II:}}
\begin{align}
   &{\rm Minimize~} Z=[\underline{z},~\overline{z}]=\sum_{i=1}^{m} \sum_{j=1}^{n}([\underline{t}_{ij},~\overline{t}_{ij}]y_{ij}+[\underline{l}_{ij},~\overline{l}_{ij}]x_{ij}) \label{0}\\
   {\rm subject~to~the~constraints;}\nonumber& \\
   &\sum_{j=1}^{n}y_{ij} \leq_{\dag}[\underline{\alpha}_{i},~ \overline{\alpha}_{i}] ~i=1,2,\ldots,{m} \label{1.2}\\
   &\sum_{i=1}^{m}y_{ij} \geq_{\dag}[\underline{\beta}_{j},~ \overline{\beta}_{j}], ~j=1,2,\ldots,{n} \label{1.3}\\
   &y_{ij}\geq{0},~i=1,2,\ldots,{m},~ j=1,2,\ldots,{n}\\
   &x_{ij} =
\begin{dcases}
       0~&{\rm if}~y_{ij}=0,\\
       1~&{\rm if}~y_{ij}>0
\end{dcases}
   \end{align}
%where the inequality relations indicated by $\leq_{\dag}$ and $\geq_{\dag}$ are defined as follows:
%\begin{align}
%w\leq_{\dag}[\alpha,\beta]\equiv\exists z\in[\alpha,\beta];~ w\leq{z},\\
%&w\geq_{\dag}[\alpha,\beta]\equiv\exists z\in[\alpha,\beta];~ w\geq{z}.
%\end{align}
   where,
   \begin{align*}
   \underline{z},~\overline{z}&: \text{left-limit and right-limit of interval $Z$, respectively } \\
   [\underline{t}_{ij},~\overline{t}_{ij}]&: \text{interval transportation cost of unit product from the $i^{th}$ source  to the $j^{th}$ destination}\\
   [\underline{l}_{ij},~\overline{l}_{ij}]&: \text{interval fixed charge between the $i^{th}$ source to the $j^{th}$ destination}\\
    [\underline{\alpha}_{i},~\overline{\alpha}_{i}]&: \text{interval availability at the $i^{th}$ source}\\
    [\underline{\beta}_{j},~\overline{\beta}_{j}]&: \text{interval demand at the $j^{th}$ destination}
   \end{align*}
\subsection{Equivalent crisp problem of IFCTP}
The objective function (\ref{0}) and constraints (\ref{1.2}-\ref{1.3}) of the IFCTP contains the interval quantities that are hard to deal with mathematically. For this reason, it is better to get a comparable crisp problem so that complex mathematical calculations can be done more easily. First of all, we convert the objective function (\ref{0}) of the IFCTP into crisp form as follows:\\
Let us consider the interval objective function (\ref{0}) of Problem-II which can be written in terms of expected value and uncertainty as; $Z=\big \langle{z_{c},z_{w}}\big \rangle$ where, $z_{c}=\big(\frac{\overline{z}+\underline{z}}{2}\big)$ and $z_{w}=\big(\frac{\overline{z}-\underline{z}}{2}\big)$.\\
In 1990, Ishibuchi and Tanaka {\cite{bakp3}} defined the center of an interval as expected value and the width of the interval as uncertainty. Since the objective function (\ref{0}) of Problem-II is the cost function which is to be minimized, so our interest is to obtain minimum expected cost with maximum precision. The left limit of the objective function (\ref{0}) can be splitted in terms of expected cost and uncertainty by definition (\ref{definition1} and \ref{definition2}) as follows:
\begin{equation}
\underline{z}=\sum_{i=1}^{m} \sum_{j=1}^{n} t_{{c}_{ij}}y_{ij}-\sum_{i=1}^{m} \sum_{j=1}^{n}t_{{w}_{ij}}y_{ij}+\sum_{i=1}^{m} \sum_{j=1}^{n} l_{{c}_{ij}}x_{ij}-\sum_{i=1}^{m} \sum_{j=1}^{n}l_{{w}_{ij}}x_{ij} \label{eq1}
\end{equation}
where, $t_{{c}_{ij}}$ and $t_{{w}_{ij}}$ are the center and width of the associated cost of $Z$ and $l_{{c}_{ij}}$ and $l_{{w}_{ij}}$ are the center and width of the fixed cost of the $Z$.\\
By equation (\ref{eq1}), we observe that minimizing (\ref{eq1}) is equivalent to minimizing the expected value (associated cost and fixed cost), but the uncertainty (associated cost and fixed cost) goes up, which is undesirable. We aim to minimize the uncertainty of the interval along with minimizing the expected value of the interval, which can be achieved by simultaneously minimizing the left limit function $(\underline{z})$ and uncertainty function $(z_{w})$. The uncertainty function can be expressed as follows:
\begin{equation}
z_{w}=\sum_{i=1}^{m} \sum_{j=1}^{n} t_{{w}_{ij}}y_{ij}+\sum_{i=1}^{m} \sum_{j=1}^{n}l_{{w}_{ij}}x_{ij} \label{eq2}
\end{equation}
Now, Using (\ref{eq1} and \ref{eq2}) we write the equivalent crisp problem of IFCTP as follows:\\ \\
{\bf{Problem-III}}
\begin{align}
&{\rm Minimize} ~\underline{z}=\sum_{i=1}^{m} \sum_{j=1}^{n} t_{{c}_{ij}}y_{ij}-\sum_{i=1}^{m} \sum_{j=1}^{n}t_{{w}_{ij}}y_{ij}+\sum_{i=1}^{m} \sum_{j=1}^{n} l_{{c}_{ij}}x_{ij}-\sum_{i=1}^{m} \sum_{j=1}^{n}l_{{w}_{ij}}x_{ij}\\
&{\rm Minimize}~z_{w}=\sum_{i=1}^{m} \sum_{j=1}^{n} t_{{w}_{ij}}y_{ij}+\sum_{i=1}^{m} \sum_{j=1}^{n}l_{{w}_{ij}}x_{ij}\\
{\rm subject~to~the~constraints;}\nonumber&\\
& \sum_{j=1}^{n}y_{ij}\leq\overline{\alpha}_{i},i=1,2,\ldots,{m} \label{eq3}\\
&\sum_{i=1}^{m}y_{ij}\geq\underline{\beta}_{j},j=1,2,\ldots,{n} \label{eq4}\\
&y_{ij}\geq{0},~i=1,2,\ldots,{m},~ j=1,2,\ldots,{n} \label{eq5}\\
&x_{ij}=
\begin{dcases}\label{eq6}
      0&~{\rm if}~y_{ij}=0,\\
      1&~{\rm if}~y_{ij}>0
   \end{dcases}.
\end{align}
\vspace{4mm}\\
The following subsection provides a solution of the equivalent bi-objective transportation problem (Problem-III). 
\subsection{Fuzzy Programming Technique for solving bi-objective transportation problem (Problem-III)}
First we find the best $L_{k}$ and worst $U_{k}$ for the $k^{th}$($k \in \{\underline{z},z_{w}\}$) objective function, where $L_{k}$ and $U_{k}$ are aspired level and highest acceptable level for the $k^{th}$ objective function, respectively. We develop a linear programming problem using membership function. The step-wise procedure of fuzzy programming technique is given as follows:
\begin{enumerate}[\bfseries{Step \arabic*:}]
\item Solve the bi-objective transportation problem (Problem-III) as a single objective problem using only one objective function at a time ignoring the other one.
\item From each solution derived in Step $1$ determine the values of both objective functions.
\item Find the best $L_{k}$ and worst $U_{k}$ for both objective corresponding to the set of solutions. Define a fuzzy membership function $\mu_{k}(z_{k})$ as follows:
  \begin{equation*}
 \mu_{k}(z_{k}) =\left\{
  \begin{array}{@{}lll@{}}
    1, & \text{if}\ z_{k}\leq L_{k} \\
    1-\frac{z_{k}-L_{k}}{U_{k}-L_{k}}, & \text{if}\ L_{k}\leq z_{k}\leq U_{k} \\
    0, & \text{if}\ z_{k}\geq L_{k}
  \end{array}\right.
\end{equation*}
Let the equivalent problem of the above vector minimum problem may be given as follows:
\begin{align*}
   &{\rm Maximize~} \lambda,\\
   {\rm subject~to~the~constraints;~}&\\
   &\lambda \leq \frac{U_{k}-z_{k}}{U_{k}-L_{k}} \\
   &{\rm and,~} (\ref{eq3}-\ref{eq6})\\
   &0 \leq {\lambda} \leq{1}
   \end{align*}
The above problem can also be expressed as follows:\\ \\
{\bf{Problem-IV:}}
\begin{align*}
   &{\rm Maximize~} \lambda,\\
   {\rm Subject~to~the~constraints;~}&\\
   &z_{k}+\lambda(U_{k}-L_{k})\leq U_{k} \\
   {\rm and,~} &(\ref{eq3}-\ref{eq6}),\\
   &0 \leq {\lambda} \leq{1}.
   \end{align*}
\item Solve Problem-IV using any method and obtain the required Pareto optimal solution.
\end{enumerate}
\subsection{Procedure for obtaining ideal solution of IFCTP (Problem-II)}
In this section we define the procedure to obtain the ideal expected value and ideal uncertainty of overall transportation cost.\\
 The objective function $Z$ of  Problem-II can be divided into two objective functions ($z_{c}$ and $z_{w}$) as follows:\\
\begin{align}
&z_{c}=\sum_{i=1}^{m} \sum_{j=1}^{n} t_{{c}_{ij}}y_{ij}+\sum_{i=1}^{m} \sum_{j=1}^{n}l_{{c}_{ij}}x_{ij},\label{eq7}\\
&z_{w}=\sum_{i=1}^{m} \sum_{j=1}^{n} t_{{w}_{ij}}y_{ij}+\sum_{i=1}^{m} \sum_{j=1}^{n}l_{{w}_{ij}}x_{ij} \label{eq8}
\end{align}
Using (\ref{eq7}) and (\ref{eq8}) construct the two linear programming problems (say Problem-V and Problem-VI) as follows:\\ \\
{\bf{Problem-V}}
\begin{align*}
&{\rm Minimize} ~z_{c}=\sum_{i=1}^{m} \sum_{j=1}^{n} t_{{c}_{ij}}y_{ij}+\sum_{i=1}^{m} \sum_{j=1}^{n}l_{{c}_{ij}}x_{ij}\\
{\rm subject~to~the~constraints;}&~\\
& (\ref{eq3}-\ref{eq6}).\nonumber
\end{align*}
{\bf{Problem-VI}}
\begin{align*}
&{\rm Minimize} ~z_{w}=\sum_{i=1}^{m} \sum_{j=1}^{n} t_{{w}_{ij}}y_{ij}+\sum_{i=1}^{m} \sum_{j=1}^{n}l_{{w}_{ij}}x_{ij}\\
{\rm subject~to~the~constraints;}&~\\
& (\ref{eq3}-\ref{eq6}).\nonumber
\end{align*}
Solve the Problem-V and Problem-VI separately and obtain their global minimums ${z_{c}}^{*}$ and ${z_{w}}^{*}$, respectively. Thus, the obtained $\big \langle{z_{c}}^{*}, {z_{w}}^{*}\big \rangle$ is the ideal solution of the Problem-II.
\vspace{4mm}\\
Now, we provide an illustration of the proposed method with an example in the following section.
\section{Numerical Illustration and comparison}
For the sake of comparison we consider the same numerical example as taken by (Safi and Razmjoo{\cite{bakp1}}). Data are given in the following Table 1.
\begin{table}[!htb]
% table caption is above the table
\caption{Cost matrix ($t_{ij},l_{ij}$)}
\centering
\begin{tabular}{l lll lll lll}
\hline\noalign{\smallskip}
&$D_{1}$ &$D_{2}$ &$D_{3}$ &$D_{4}$ &Supply\\
\noalign{\smallskip}\hline\noalign{\smallskip}
$S_{1}$ &([4,8],[10,30]) &([8,12],[19,25]) &([9,11],[19,25]) &([8,10],[20,30]) &[30,33] \\
\noalign{\smallskip}
$S_{2}$ &([10,18],[16,20]) &([10,12],[15,25]) &([11,15],[25,55]) &([5,7],[38,40]) &[27,28] \\
\noalign{\smallskip}
$S_{3}$ &([7,19],[10,20]) &([8,12],[22,30]) &([8,14],[30,50]) &([13,17],[20,22]) &[22,25] \\
\noalign{\smallskip}\noalign{\smallskip}
Demand &[20,21] &[19,24] &[23,24] &[20,22] \\
\hline\noalign{\smallskip}\noalign{\smallskip}
\end{tabular}
\end{table}\\
First, we have converted the given problem into an equivalent crisp problem using the approach described in section 3.\\ \\
{\bf{Problem-VII}}
\begin{align*}
&{\rm Minimize}~\underline{z}=\sum_{i=1}^{3} \sum_{j=1}^{4} t_{{c}_{ij}}y_{ij}-\sum_{i=1}^{3} \sum_{j=1}^{4}t_{{w}_{ij}}y_{ij}+\sum_{i=1}^{3} \sum_{j=1}^{4} l_{{c}_{ij}}x_{ij}-\sum_{i=1}^{3} \sum_{j=1}^{4}l_{{w}_{ij}}x_{ij}\\
{\rm{equivalently}},\nonumber\\
&{\rm Minimize}~\underline{z}=\sum_{i=1}^{3} \sum_{j=1}^{4} \underline{{t}}_{ij}y_{ij}+\sum_{i=1}^{3} \sum_{j=1}^{4} \underline{{l}}_{ij}x_{ij},\\
&{\rm and~ minimize}~z_{w}=\sum_{i=1}^{3} \sum_{j=1}^{4} t_{{w}_{ij}}y_{ij}+\sum_{i=1}^{m} \sum_{j=1}^{n}l_{{w}_{ij}}x_{ij}\\
&{\rm subject~to~the~constraints;}\nonumber\\
& \sum_{j=1}^{4}y_{1j}\leq33, \sum_{j=1}^{4}y_{2j}\leq28, \sum_{j=1}^{4}y_{3j}\leq25; \\
&\sum_{i=1}^{3}y_{i1}\geq20, \sum_{i=1}^{3}y_{i2}\geq19, \sum_{i=1}^{3}y_{i3}\geq23, \sum_{i=1}^{3}y_{i4}\geq20; \\
&y_{ij}\geq{0},~i=1,2,3,~ j=1,2,3,4\\
&x_{ij}=
\begin{dcases}
      0&~{\rm if}~y_{ij}=0,\\
      1&~{\rm if}~y_{ij}>0
   \end{dcases}
\end{align*}
where,$$
\ \underline{{t}}_{ij}=
\begin{bmatrix}
4& 8 &9 &8\\
10 &10 &11 &5\\
7 &8 &8 &13\\
\end{bmatrix},
\ \underline{{l}}_{ij}=
\begin{bmatrix}
10& 19 &19 &20\\
16 &15 &25 &38\\
10 &22 &30 &20\\
\end{bmatrix}
$$
$$
\ t_{{w}_{ij}}=
\begin{bmatrix}
2 &2 &1 &1\\
4 &1 &2 &1\\
6 &2 &3 &2\\
\end{bmatrix},
\ l_{{w}_{ij}}=
\begin{bmatrix}
10& 3 &3 &5\\
2 &5 &15 &1\\
5 &4 &10 &1\\
\end{bmatrix}.
$$
The above problem (Problem-VII) has been coded in LINGO 19.0 and solved. The LINGO codes are given in the following frame:
\begin{Verbatim}[frame=single]
sets:
 supply /1..3/:ar;
 demand /1..4/:bl;
 link(supply, demand): tl,ll,tw,lw,x,y;
endsets
data:
ar = 33 28 25;
bl = 20 19 23 20;
tl  = 4 8 9 8 10 10 11 5 7 8 8 13;
ll  = 10 19 19 20 16 15 25 38 10 22 30 20;
tw  = 2 2 1 1 4 1 2 1 6 2 3 2;
lw  =10 3 3 5 2 5 15 1 5 4 10 1;
U1=787; L1=640;
U2=190; L2=163;
enddata
max=LAMBDA;
@sum(supply(i):@sum(demand(j):tl(i,j)*y(i,j)+ll(i,j)*x(i,j)))+LAMBDA*(U1-L1)<=U1;
@sum(supply(i):@sum(demand(j):tw(i,j)*y(i,j)+lw(i,j)*x(i,j)))+LAMBDA*(U2-L2)<=U2;
@for(demand(j):@sum(supply(i):y(i,j))>=bl(j););
@for(supply(i):@sum(demand(j):y(i,j))<=ar(i););
@for(link(i,j):x(i,j)=@if(y(i,j)#gt#0,1,0));
LAMBDA>=0;
END
\end{Verbatim}
The Pareto optimal solution of the problem is $y_{11}=20$, $y_{13}=13$, $y_{22}=4.95$, $y_{24}=20$, $y_{32}=14.04$, $y_{33}=10$ and $x_{11}=1$, $x_{13}=1$, $x_{22}=1$, $x_{24}=1$, $x_{32}=1$, $x_{33}=1$ and remaining $y_{ij}$'s and $x_{ij}$'s are equal to zero; $Z$=$[672.82,1010.88]$=
${\big \langle{841.85,169.03}\big \rangle}$. \\
We obtain the ideal expected value and ideal uncertainty of the given problem as ${\big \langle{830,163}\big \rangle}$ using the procedure given in section 3.4.\\Using definition ({\ref{def1}}) the distance from ${\big \langle{830,163}\big \rangle}$ to ${\big \langle{841.85,169.03}\big \rangle}$ is $13.29$ .\\
\subsection{Comparison}
This paper discusses a solution approach of an IFCTP in which the involved interval numbers represent uncertain parameters. A new approach for solving such IFCTP has been developed and compared with an existing approach suggested by Safi and Razmjoo \cite{bakp1}. Two existing numerical examples have been taken from Safi and Razmjoo \cite{bakp1} and solved using the proposed approach. A detailed solution has been provided for one of the two numerical examples, along with its LINGO codes, to demonstrate the proposed method. Table 2 clearly shows that the results provided by the proposed method are closer to their respective ideal solutions than that of Safi and Razmjoo \cite{bakp1} ( in both the numerical examples). Thus, the decision maker has the opportunity to use the proposed method in order to get better results.
\FloatBarrier
\begin{table}[h!]
 \centering
\caption{\bf{Comparison of proposed method with the method given by Safi and Razmjoo {\cite{bakp1}}}}
\centering
  \label{table1}
\begin{tabular}{l c c c c l }
 \hline
 & & \multicolumn{2}{c}{Solution} &\multicolumn{2}{m{3.5cm}} {Distance of the solution from ideal solution} \\
 \hline
 Problem & \multicolumn{1}{m{1.5cm}} {Ideal Solution $\big \langle{z_{c}^{*},{z_{w}^{*}}}\big \rangle$}  & \multicolumn{1}{m{1cm}}{Safi and Razmjoo {\cite{bakp1}}}&\multicolumn{1}{m{1.5cm}}{Proposed Method}  & \multicolumn{1}{m{1.5cm}}{Safi and Razmjoo{\cite{bakp1}}} &\multicolumn{1}{m{1.5cm}}{Proposed Method} \\
 \hline\noalign{\smallskip}\noalign{\smallskip}
 1 & \big \langle{830,163}\big \rangle  & \multicolumn{1}{m{2.2cm}} {[640,1020] \big \langle{830,190}\big \rangle}& \multicolumn{1}{m{2.2cm}} {[672.82,1010.88] \big \langle{841.85,169.03}\big \rangle}  & 27 & \bf{13.29} \\
 \noalign{\smallskip}\noalign{\smallskip}
 2 & \big \langle{734.5,26.75}\big \rangle  & \multicolumn{1}{m{2.2cm}} {[734,770] \big \langle{752,18}\big \rangle} & \multicolumn{1}{m{2.2cm}} {[712.5,767.5] \big \langle{740,27.5}\big \rangle}  & 8.76 & \bf{5.55}\\
 \hline\noalign{\smallskip}\noalign{\smallskip}
 \end{tabular}
 \end{table}
 \FloatBarrier
 \section{Concluding remarks}
 Transportation problems involving interval-valued parameters hold significant importance in various real-world scenarios. These problems address the inherent uncertainty that often surrounds transportation and logistics decisions. Unlike traditional deterministic models, interval-valued parameters allow for the representation of ranges rather than fixed values, accommodating the variability and imprecision inherent in factors such as travel times, costs, and demand. This approach enables decision-makers to devise more robust and resilient transportation strategies, considering a range of potential outcomes rather than relying on single-point estimates. By embracing uncertainty and incorporating interval-valued parameters, transportation planners can create solutions that are better equipped to handle unforeseen disruptions, market fluctuations, and other dynamic challenges, ultimately leading to improved decision-making, resource allocation, and overall system efficiency in the complex realm of transportation management. \par 
 
 This article investigates a solution strategy for the Interval Fixed Charge Transportation Problem (IFCTP), in which interval numbers represent uncertain parameters. A novel method for addressing this IFCTP has been developed and compared to one proposed by Safi and Razmjoo {\cite{bakp1}}. The study demonstrates the effectiveness of the newly developed method by applying it to two numerical examples initially presented by Safi and Razmjoo {\cite{bakp1}}.

\end{document}